\documentclass{gtart_h}

\def\ifplaintex{\expandafter\ifx\csname documentclass\endcsname\relax}

\def\gtp{{\mathsurround=0pt\it $\cal G\mskip-2mu$eometry \&\ 
$\cal T\!\!$opology $\cal P\!$ublications}}  % GT publications

\def\recd{{\small Received:\qua\receiveddate\ifx\reviseddate\relax
\else\qquad Revised:\qua\reviseddate\fi\par}} 

\def\lognumber#1{\def\thelognumber{#1}}
\def\volumenumber#1{\def\thevolumenumber{#1}}
\def\volumeyear#1{\def\thevolumeyear{#1}}
\def\papernumber#1{\def\thepapernumber{#1}}
\def\pagenumbers#1#2{\def\startpage{#1}\def\finishpage{#2}}
\def\published#1{\def\publishdate{#1}}

\def\received#1{\def\receiveddate{#1}}
\def\revised#1{\def\reviseddate{#1}}
\def\accepted#1{\def\accepteddate{#1}}

\long\def\asciiabstract#1{\long\def\theasciiabstract{#1}}

\let\\\par\let\thelognumber\relax\let\thevolumenumber\relax
\let\thepapernumber\relax\let\thevolumeyear\relax\let\startpage\relax
\let\finishpage\relax\let\publishdate\relax\let\receiveddate\relax
\let\reviseddate\relax\let\accepteddate\relax\let\theasciititle\relax
\let\theasciiauthors\relax
\let\theasciiabstract\relax

\let\theasciiemail\relax

\ifplaintex
\font\logobig=cmssbx10 scaled 3836
\font\logomed=cmssbx10 scaled 2557
\else
\font\logobig=cmssbx10 scaled 4200
\font\logomed=cmssbx10 scaled 2800
\fi

\long\def\makeagttitle{   %%% start of definition of \makeagttitle
\count0=\startpage
\agt\hfill      %   Journal title (top left) 
\hbox to 45truept{\vbox to 0pt{\vglue -13truept{\logomed A\kern -.37em{\logobig 
T}\kern -.38em G}\vss}\hss}
\break
{\small Volume \thevolumenumber\ (\thevolumeyear)
\startpage--\finishpage\nl
Published: \publishdate}

\vglue .25truein

{\parskip=0pt\leftskip 0pt plus
1fil\def\\{\par\smallskip}{\Large\bf\thetitle}\par\medskip} \vglue
0.05truein

{\parskip=0pt\leftskip 0pt plus 1fil\def\\{\par}{\sc\theauthors}
\par\medskip}%
 
\vglue 0.03truein 

{\small\leftskip 25truept\rightskip 25truept{\bf Abstract}\stdspace\theabstract

{\bf AMS Classification}\stdspace\theprimaryclass
\ifx\thesecondaryclass\relax\else; \thesecondaryclass\fi\par
{\bf Keywords}\stdspace \thekeywords\par}\vglue 7truept

}   %%%% end of definition of \makeagttitle

\ifplaintex
\font\phead=cmsl9 scaled 950
\font\pnum=cmbx10 scaled 913
\font\pfoot=cmsl9 scaled 950
\headline{\vbox to 0pt{\vskip -4.5mm\line{\small\phead\ifnum
\count0=\startpage ISSN 1472-2739 (on-line) 1472-2747 (printed)
\hfill {\pnum\folio}\else\ifodd\count0\def\\{ }% 
\ifx\theshorttitle\relax\thetitle\else\theshorttitle\fi\hfill{\pnum\folio}
\else\def\\{ and }{\pnum\folio}\hfill\ifx\theshortauthors\relax\theauthors
\else\theshortauthors\fi\fi\fi}\vss}}
\footline{\vbox to 0pt{\vglue 0mm\line{\small\pfoot\ifnum\count0=\startpage
\copyright\ \gtp\hfill\else
\agt, Volume \thevolumenumber\ (\thevolumeyear)\hfill\fi}\vss}}
\else
\font=cmsl9 scaled 1050
\font\lnum=cmbx10 
\font=cmsl9 scaled 1050
\makeatletter
\def\@oddhead{{\small\lhead\ifnum\count0=\startpage ISSN 1472-2739 
(on-line) 1472-2747 (printed)\hfill {\lnum\number\count0}\else\ifodd\count0
\def\\{ }\ifx\theshorttitle\relax \thetitle \else\theshorttitle\fi\hfill
{\lnum\number\count0}\else\def\\{ and }{\lnum\number\count0}
\hfill\ifx\theshortauthors\relax 
\theauthors\else\theshortauthors\fi\fi\fi}}\def\@evenhead{\@oddhead}
\def\@oddfoot{\small\lfoot\ifnum\count0=\startpage\copyright\ \gtp\hfill\else
\agt, Volume \thevolumenumber\ (\thevolumeyear)\hfill\fi}
\def\@evenfoot{\@oddfoot}
\makeatother
\fi
\let\maketitlepage\makeagttitle
\let\makeshorttitle\maketitlepage
\let\maketitle\maketitlepage

\newwrite\gtoutfile
\long\gdef\makeheadfile{  %%% start of definition of \makeheadfile
{\def\\{, }\def\s{ }
\immediate\openout\gtoutfile head.xxx
\immediate\write\gtoutfile{Proxy-for: \ifx\theasciiauthors\relax
\theauthors\else\theasciiauthors\fi\s<\ifx\theasciiemail\relax\theemail\else\theasciiemail\fi>}
\immediate\write\gtoutfile{\noexpand\\}
\immediate\write\gtoutfile{Authors: \ifx\theasciiauthors\relax
\theauthors\else\theasciiauthors\fi}
{\def\\{ }\immediate\write\gtoutfile{Title: \ifx\theasciititle\relax
\thetitle\else\theasciititle\fi}}
\immediate\write\gtoutfile{Subj-class: GT or SG, GR etc}
\immediate\write\gtoutfile{MSC-class: \theprimaryclass\ifx\thesecondaryclass\relax\else, \thesecondaryclass\fi}
\immediate\write\gtoutfile{Journal-ref: Algebr. Geom. Topol. \thevolumenumber\s
(\thevolumeyear) \startpage-\finishpage}
\immediate\write\gtoutfile{Comments: Published by Algebraic and
Geometric Topology at}
\immediate\write\gtoutfile{\s\s\s  http://www.maths.warwick.ac.uk/agt/AGTVol\thevolumenumber/agt-\thevolumenumber-\thepapernumber.abs.html}
\immediate\write\gtoutfile{\noexpand\\}
\immediate\write\gtoutfile{}
\ifx\theasciiabstract\relax
\immediate\write\gtoutfile{\theabstract}\else
\immediate\write\gtoutfile{\theasciiabstract}\fi
\immediate\write\gtoutfile{}
\immediate\write\gtoutfile{\noexpand\\}
\immediate\write\gtoutfile{}
\immediate\closeout\gtoutfile}}  %%% end of definition of \makeheadfile

\def\maketitlepage{\makeagttitle\makeheadfile}
\let\makeshorttitle\maketitlepage
\let\maketitle\maketitlepage

\lognumber{28}
\volumenumber{4}
\volumeyear{2004}
\papernumber{28}
\received{19 September 2003}
\pagenumbers{603}{622}
\published{20 August 2004}
\revised{14 May 2004}
\accepted{2 August 2004}

\usepackage{amsmath,amssymb}

\numberwithin{equation}{section}
 
\theoremstyle{plain}

\newtheorem{Th}[equation]{Theorem}
\newtheorem{Prop}[equation]{Proposition}
\newtheorem{Le}[equation]{Lemma}
\newtheorem{Cor}[equation]{Corollary}

\theoremstyle{remark}

\theoremstyle{definition}

\newcommand{\Z}{\mathbb Z}
\newcommand{\R}{\mathbb R}
\newcommand{\FCC}{{\rm FCC} }
\newcommand{\Link}{{\rm Link}}
\newcommand{\image}{{\rm image}}
\newcommand{\CAT}{{\rm CAT}}
\newcommand{\FCCs}{{\rm FCC}s }

\def\ra{\rightarrow}

\def\i{\infty}

\begin{document}
\title{Foldable cubical complexes of\\nonpositive curvature}
\author{Xiangdong Xie}
\address{Dept.\ of Math.\ Sciences, University of 
Cincinnati, Cincinnati, OH 45221, USA}
\email{xiexg@ucmail.uc.edu}

\begin{abstract}
We study finite foldable cubical complexes of nonpositive curvature
(in the sense of A.D.\ Alexandrov).  We show that such a complex $X$
admits a graph of spaces decomposition.  It is also shown that when
$\dim X=3$, $X$ contains a closed rank one geodesic in the
$1$-skeleton unless the universal cover of $X$ is isometric to the
product of two $\CAT(0)$ cubical complexes.  
\end{abstract}

\asciiabstract{%
We study finite foldable cubical complexes of nonpositive curvature
(in the sense of A.D. Alexandrov).  We show that such a complex X
admits a graph of spaces decomposition.  It is also shown that when
dim X=3, X contains a closed rank one geodesic in the 1-skeleton
unless the universal cover of X is isometric to the product of two
CAT(0) cubical complexes.}


\maketitle

\section{Introduction}

We study finite foldable cubical complexes with nonpositive curvature
(in the sense of A.D.\ Alexandrov), including the rank rigidity
problem of such complexes.  Foldable cubical complexes have been
studied by W. Ballmann and J. Swiatkowski in \cite{BSw}. Our notion of
foldable cubical complexes is slightly more general than that of
\cite{BSw} since we do not require gallery connectedness.  D. Wise has
also studied a class of 2-dimensional cubical complexes
($\mathcal{V}H$-complexes in \cite{W}) which are closely related to
foldable cubical complexes.

A \emph{cubical complex} $X$ is a CW-complex formed by gluing unit
 Euclidean cubes together along faces via isometries.  We require that
 all the cubes inject into $X$ and the intersection of the images of
 two cubes is either empty or equals the image of a cube.  The image
 of a $k$-dimensional unit Euclidean cube in $X$ is called a $k$-cube,
 and a $1$-cube is also called an edge.  Let $d$ be the path
 pseudometric on $X$.  When $X$ is finite dimensional, $d$ is actually
 a metric and turns $X$ into a complete geodesic space, see \cite{Br}.

Let $X$ be a cubical complex.  $X$ is called \emph{dimensionally
    homogeneous} if there is an integer $n\ge 1$ such that each cube
    of $X$ is a face of some $n$-cube. A \emph{folding} of $X$ is a
    combinatorial map $f: X\ra C$ onto an $n$-cube $C$ such that the
    restriction of $f$ on each cube is injective.  $X$ is
    \emph{foldable} if it admits a folding.  $X$ has \emph{nonpositive
    curvature} in the sense of A.D.\ Alexandrov if and only if all the
    vertex links are flag complexes \cite{BH}.  Recall that a
    simplicial complex is a flag complex if any finite set of vertices
    that is pairwise joined by edges spans a simplex.  By a
    \textsl{FCC} we mean a connected foldable cubical complex that is
    dimensionally homogeneous, geodesically complete (definition given
    below) and has nonpositive curvature.

 Our first observation about FCCs is the following:

\medskip

\noindent
{\bf{Proposition \ref{graphdeco}}}\qua \textsl{Let $X$ be a \FCC of
dimension $n$. Then $X$ admits the structure of a graph of spaces,
where all the vertex and edge spaces are $(n-1)$-dimensional \FCCs and
the maps from edge spaces to vertex spaces are combinatorial
immersions.}

\medskip

Notice Proposition \ref{graphdeco} offers the potential for proving
 statements about FCCs by inducting on dimension.

We shall give two applications of the above observation.  The first
 concerns the rank rigidity problem for $\CAT(0)$ spaces.  Let $X$ be
 a metric space with nonpositive curvature. A curve $\sigma: I\ra X$
 is a \emph{geodesic} if it has constant speed and is locally
 distance-minimizing.  We say $X$ is \emph{geodesically complete} if
 every geodesic $\sigma: I\ra X$ can be extended to a geodesic $\tilde
 \sigma:\R\ra X$.  A \emph{$\CAT(0)$ space} is a simply connected
 complete geodesic space with nonpositive curvature.  Let $Y$ be a
 geodesically complete $\CAT(0)$ space.  We say $Y$ has higher rank if
 each geodesic $\sigma: \R\ra Y$ is contained in a flat plane (the
 image of an isometric embedding ${\R}^2\ra Y$).  Otherwise we say $Y$
 has rank one.  There are two main classes of higher rank $\CAT(0)$
 spaces: symmetric spaces and Euclidean buildings.  The following
 conjecture is still open (\cite{BBu}, \cite{BBr2}):

\medskip

\noindent
{\bf{Rank Rigidity Conjecture}}\qua \textsl{Let $Y$ be a geodesically
complete $\CAT(0)$ space and $\Gamma$ a group of isometries of $Y$
whose limit set is the entire ideal boundary of $Y$.}

(1)\qua \textsl{If $Y$ has higher rank, then $Y$ isometrically splits
 unless $Y$ is a symmetric space or an Euclidean building;}

(2)\qua \textsl{If $Y$ has rank one, then $Y$ contains a periodic rank
one geodesic.}

\medskip

Recall a complete geodesic $\sigma: \R\ra Y$ is a \emph{periodic rank
     one geodesic} if $\sigma$ does not bound a flat half-plane and
     there is some $\gamma\in \Gamma$ and $c>0$ such that
     $\gamma(\sigma(t))=\sigma(t+c)$ for all $t\in \R$.  The
     conjecture holds if the action of $\Gamma$ on $Y$ is proper and
     cocompact and $Y$ is a Hadamard manifold \cite{B} or a
     $2$-dimensional polyhedron with a piecewise smooth metric
     \cite{BBr}.  Claim (1) of the conjecture also holds when $Y$ is
     a $3$-dimensional piecewise Euclidean polyhedron and $\Gamma$
     acts on $Y$ properly and cocompactly \cite{BBr2}.

\medskip

\noindent
{\bf{Theorem \ref{rigidity}}}\qua \textsl{Let $X$ be a finite \FCC of
dimension $3$ with universal cover $\tilde X$ and group of deck
transformations $\Gamma$.}

(1)\qua \textsl{If $\tilde X$ has higher rank, then $\tilde X$ is
 isometric to the product of two $\CAT(0)$ \FCCs.}

(2)\qua \textsl{If $\tilde X$ has rank one, then $\tilde X$ contains a
  periodic rank one geodesic in the $1$-skeleton.}

\medskip

As the second application, we address the Tits alternative question
  for finite \FCCs.  The result has been established by Ballmann and
  Swiatkowski \cite{BSw}.  We give a new and very short proof.

\medskip

\noindent
{\bf{Theorem \ref{tits}}}\qua \textsl{Let $X$ be a finite \FCC.  Then
any subgroup of $\pi_1(X)$ either contains a free group of rank two or
is virtually free abelian.}

\medskip

The paper is organized as follows.  In Section \ref{fold} we recall
 some basic facts about FCCs and show that a FCC has a graph of spaces
 structure.  In Section \ref{existence} the rank rigidity problem for
 $3$-dimensional FCCs is discussed.  In Section \ref{titsal} we give a
 new proof of the Tits alternative for FCCs.

The author would like to thank the referee and the editor for many
helpful comments.

\section{Foldable cubical complexes}\label{fold}

\subsection{Locally convex subcomplexes} \label {subcom}

In this section we show that a FCC has many locally convex
subcomplexes.  The results in this section are more or less known (see
\cite{BSw}, \cite{C}, \cite{DJS}).  We include the proofs only for
completeness.

 For any locally finite piecewise Euclidean complex $Y$ and any $y\in
Y$, the link $\Link(Y, y)$ is piecewise spherical.  We let $d_y$ be the
induced path metric on $\Link(Y, y)$.  If $Y$ has nonpositive
curvature, then $\Link(Y, y)$ is a $\CAT(1)$ space.  Let $Y$ be a
locally finite piecewise Euclidean complex with nonpositive curvature
and $Z\subset Y$ a subcomplex.  For $z\in Z$, a subset $L(Z,z)\subset
\Link(Y,z)$ is defined as follows: a point $\xi\in \Link(Y,z)$ belongs
to $L(Z,z)$ if there is a nontrivial geodesic segment $zx$ of $Y$ with
$zx\subset Z$ such that the initial direction of $zx$ at $z$ is $\xi$.

For any metric space $Z$, the \emph{Euclidean cone} over $Z$ is the
metric space $C(Z)$ defined as follows.  As a set $C(Z)=Z\times
[0,\i)/{Z\times \{0\}}$.  The image of $(z,t)$ is denoted by $tz$.
$d(t_1z_1, t_2z_2)=t_1+t_2$ if $d(z_1, z_2)\ge \pi$, and $d(t_1z_1,
t_2z_2)=\sqrt{t_1^2+t_2^2-2t_1t_2\cos(d(z_1,z_2))}$ if $d(z_1, z_2)\le
\pi$.  The point $O:=Z\times \{0\}$ is called the cone point of
$C(Z)$.

   Recall a subset $A\subset M $ of a $\CAT(1)$ metric space $M$ is
  $\pi$-convex if for any $a, b\in A$ with $d(a,b)<\pi$ the geodesic
  segment $ab$ lies in $A$.
 
\begin{Le}\label{cone}
{Let $Y$ be a locally finite piecewise Euclidean complex with
 nonpositive curvature and $Z\subset Y$ a subcomplex.  Then $Z$ is
 locally convex in $Y$ if and only if for each $z\in Z$, $L(Z,z)$ is
 $\pi$-convex in $\Link(Y,z)$.}
\end{Le}

\begin{proof}
For any $y\in Y$ let $C(\Link(Y, y))$ be the Euclidean cone over
$\Link(Y,y)$ and $O$ the cone point.  For any $r>0$ let $\bar
B(y,r)\subset Y$ and $\bar B(O, r)\subset C(\Link(Y, y))$ be the closed
metric balls with radius $r$.  For any subset $A\subset \Link(Y, y)$,
let $C_r(A)\subset C(\Link(Y, y))$ be the subset consisting of points
of the form $ta$, $t\le r$ and $a\in A$.

Since $Y$ is a locally finite piecewise Euclidean complex and
$Z\subset Y$ is a subcomplex, for each $z\in Z$ there is some $r>0$
and an isometry $g: \bar B(z,r)\ra \bar B(O,r)$ such that $g(\bar
B(z,r)\cap Z)=C_r(L(Z,z))$.  Now it is easy to see that $ \bar
B(z,r)\cap Z$ is convex in $\bar B(z,r)$ if and only if $L(Z,z)$ is
$\pi$-convex in $\Link(Y,z)$.
\end{proof}

Let $X$ be a FCC and $f: X\ra C$ a fixed folding.  Two edges $e_1$ and
 $e_2$ of $X$ are equivalent if $f(e_1)$ and $f(e_2)$ are parallel in
 $C$.  Let $E_1,\cdots, E_n$ be the equivalence classes of the edges
 of $X$.  For each nonempty subset $T\subset \{1, 2,\cdots, n\}$ we
 define a subcomplex $X_T$ of $X$ as follows: a $k$-cube belongs to
 $X_T$ if all its edges belong to $E_T:=\cup_{i\in T}E_i$.  $X_T$ is
 disconnected in general.  We shall see that each component of $X_T$
 is locally convex in $X$.

To prove that the components of $X_T$ are locally convex in $X$ we
 also need the following lemma.  Recall a spherical simplex is all
 right if all its edges have length $\pi/2$, and a piecewise spherical
 simplicial complex $K$ is all right if all its simplices are all
 right.  When $K$ is finite dimensional, the path pseudometric on $K$
 is a metric that makes $K$ a complete geodesic space \cite{Br}.

\begin{Le}[\cite{BH}, p.211]\label{pi}
{Let $K$ be a finite dimensional all right spherical complex and $v\in
K$ a vertex. If $\sigma: [a,b]\ra K$ is a geodesic such that
$\sigma(a), \sigma(b)\notin B(v,\pi/2)$, then each component of
$B(v,\pi/2)\cap \sigma$ has length $\pi$.  }
\end{Le}

A dimensionally homogeneous cubical complex $Z$ of dimension $n$
\emph{has no boundary} if each $(n-1)$-cube is contained in the
boundaries of at least two \linebreak $n$-cubes.  Similarly, a
dimensionally homogeneous simplicial complex $Z$ of \linebreak
dimension $n$ has no boundary if each $(n-1)$-simplex is contained in
the boundaries of at least two $n$-simplices.

The following proposition is a consequence of Proposition 1.7.1 of
\cite{DJS}.  It also follows from Lemmas 1.1 and 1.3 of \cite{C}.  In
addition, W. Ballmann and J. Swiatkowski have made the same
observation(\cite{BSw}, Lemma 3.2(4) and the paragraph at the bottom
of p.69 and the top of p.70).

\begin{Prop}\label{convex}
{Let $X$ be a \FCC and $f: X\ra C$ a fixed folding onto an $n$-cube.
   Then for any nonempty $T\subset \{1,2 \cdots, n\}$, each component
   of $X_T$ is locally convex in $X$. Furthermore, each component of
   $X_T$ is also a \FCC.  }
\end{Prop}
\begin{proof}
Let $Z$ be a component of $X_T$. For any $z\in Z$ we need to show that
$L(Z,z)\subset \Link(X, z)$ is $\pi$-convex. First assume $z$ is a
vertex.  Then $\Link(X,z)$ is an all right flag complex. By the
definition of $X_T$ we see $L(Z,z)$ is a full subcomplex of
$\Link(X,z)$, that is, a simplex of $\Link(X,z)$ lies in $L(Z,z)$ if and
only if all its vertices lie in $L(Z,z)$.  Let $\xi, \eta\in L(Z,z)$
with $d_z(\xi, \eta)<\pi$. Assume $\xi\eta \not\subseteq L(Z,z)$. Then
there is some $\xi'\in \xi\eta-L(Z,z)$.  Let $\Delta$ be the smallest
simplex of $\Link(X,z)$ containing $\xi'$ and $\omega_1, \cdots,
\omega_k$ be its vertices.  Then $\xi'\in B(\omega_j,\pi/2)$ for all
$1\le j\le k$.  Since $L(Z,z)$ is a full subcomplex and $\xi'\notin
L(Z,z)$, there exists some $j$, $1\le j\le k$, with $\omega_j\notin
L(Z,z)$.  We may assume $\omega_1\notin L(Z,z)$.  By the definition of
$X_T$, $\omega_1$ corresponds to some edge $e\in E_i$ with $i\notin
T$.  It follows that $\xi, \eta\notin B(\omega_1, \pi/2)$.  Then Lemma
\ref{pi} implies that $\xi\eta\cap B(\omega_1, \pi/2)$ has length
$\pi$, contradicting to the fact that $d_z(\xi, \eta)<\pi$.

When $z$ is not a vertex, there is an obvious all right simplicial
complex structure on $\Link(X, z)$ where the vertices in $\Link(X,z)$
are represented by geodesic segments parallel to the edges in $X$.  In
this case $L(Z,z)$ is still a full subcomplex of $\Link(X,z)$ and the
above argument applies.

It is clear that $Z$ is a foldable cubical complex that is
 dimensionally homogeneous.  Since $X$ is geodesically complete, it
 has no boundary.  It follows that $Z$ also has no boundary.  $Z$ is
 locally convex in $X$ implies $Z$ has nonpositive curvature in the
 path metric. By Proposition \ref{gcomplete} $Z$ is geodesically
 complete and therefore is a FCC.
\end{proof}

\begin{Prop}\label{gcomplete}
{Let $Z$ be either a piecewise Euclidean complex that has nonpositive
curvature, or a piecewise spherical complex that is locally $\CAT(1)$.
Assume $Z$ is locally finite and dimensionally homogeneous.  Then $Z$
is geodesically complete if and only if it has no boundary.  }
\end{Prop}
\begin{proof}
 It is clear that if $Z$ is geodesically complete then it has no
 boundary.  To prove the other direction, we assume $Z$ has no
 boundary.  We first notice that $Z$ is geodesically complete if and
 only if for any $z\in Z$ and any $\xi\in \Link(Z,z)$ there is some
 $\eta\in \Link(Z,z)$ with $d_z(\xi,\eta)\ge \pi$.  The link
 $\Link(Z,z)$ is a finite piecewise spherical complex that is
 $\CAT(1)$, dimensionally homogeneous and has no boundary.  We prove
 the following statement by induction on the dimension: Let $Y$ be an
 $n$-dimensional finite piecewise spherical complex that is $\CAT(1)$,
 dimensionally homogeneous and has no boundary.  Then for any $x\in
 Y$, there is some $y\in Y$ with $d(x,y)\ge \pi$.

When $n=1$, $Y$ is a finite graph and the claim is clear.  Let $n=\dim
 Y\ge 2$ and suppose there is a point $x\in Y$ such that $d(x,y)<\pi$
 for all $y\in Y$.  Since $Y$ is finite, there is some $y_0\in Y$ with
 $d(x,y)\le d(x,y_0)$ for all $y\in Y$.  As $Y$ is $\CAT(1)$ and
 $d(x,y_0)<\pi$, there is a unique minimizing geodesic
 $\sigma:[0,d(x,y_0)]\ra Y$ from $y_0$ to $x$. Let $\xi\in \Link(Y,
 y_0)$ be the point represented by $\sigma$.  Now $\Link(Y, y_0)$ has
 dimension $n-1$ and the induction hypothesis implies that there is
 some $\eta\in \Link(Y, y_0)$ with $d_{y_0}(\xi, \eta)\ge \pi$, where
 $d_{y_0}$ is the path metric on $\Link(Y, y_0)$.  Hence $\sigma$ can
 be extended to a geodesic $\tilde \sigma: [-\epsilon, d(x,y_0)]\ra Y$
 that contains $y_0$ in the interior.  As $Y$ is $\CAT(1)$, $\tilde
 \sigma$ is minimizing for small enough $\epsilon$, contradicting to
 the choice of $y_0$.
\end{proof}

\subsection{Graph of spaces decomposition} \label{decom}

In this section we show that a FCC admits decomposition as a graph of
spaces, as defined in \cite{SW}.

Let $X$ be a FCC and $f:X\ra C_0$ a fixed folding onto an
$n$-cube. Then the set $E$ of edges of $X$ is a disjoint union
$E=\amalg_{i=1}^{n}E_i$, see Section \ref{subcom}.  For each $i$ with
$1\le i\le n$, let $T_i=\{1,2,\cdots, n\}-\{i\}$.  Then the components
of $X_{T_i}$ are $(n-1)$-dimensional FCCs and are locally convex in
$X$.

For each $n$-cube $C$ of $X$, let $e\in E_i$ be an edge of $C$ and
 $C_i\subset C$ the hyperplane in $C$ containing the midpoint of $e$
 and perpendicular to $e$.  It is clear that $C_i\subset C$ does not
 depend on $e$ and is isometric to a $(n-1)$-dimensional unit
 Euclidean cube.  Set $H_i=\cup C_i$, where $C$ varies over all
 $n$-cubes of $X$.  $H_i$ is not connected in general.  An argument
 similar to the one in Section \ref{subcom} shows that each component
 of $H_i$ is locally convex in $X$.  $H_i$ has a natural FCC
 structure, where each $C_i$ is a $(n-1)$-cube.

It is not hard to see that each component of $X-X_{T_i}$ is isometric
to $Y\times (0,1)$, where $Y$ is a component of $H_i$.  Let $\{Y_1,
\cdots, Y_k\}$ be the set of components of $H_i$. Then $X$ can be
obtained from $X_{T_i}$ by attaching $Y_j\times [0,1]$ along
$Y_j\times \{0\}$ and $Y_j\times \{1\}$, $1\le j\le k$.  That is, $X$
has the structure of a graph of spaces \cite{SW} for each $i$, $1\le
i\le n=\dim X$.  Now we make it more precise.  The base graph $G_i$ is
as follows.  The vertex set $\{v_B\}$ of $G_i$ is in one-to-one
correspondence with the set $\{B\}$ of components of $X_{T_i}$, and
the edge set $\{e_Y\}$ is in one-to-one correspondence with the set
$\{Y\}$ of components of $H_i$.  For each edge $e_Y$, consider the
component $Y\times (0,1)$ of $X-X_{T_i}$ corresponding to $Y$.  The
closure of $Y\times (0,1)$ in $X$ has nonempty intersection with one
or two components of $X_{T_i}$. Let $B_0, B_1$ be these components of
$X_{T_i}$ (we may have $B_0=B_1$).  Then the edge $e_Y$ connects the
vertices $v_{B_0}$ and $v_{B_1}$.  The vertex space corresponding to
$v_B$ is the component $B$ of $X_{T_i}$ and the edge space
corresponding to $e_Y$ is the component $Y$ of $H_i$.

We notice that the base graph $G_i$ is connected: Let $v_B$, $v_{B'}$
be two vertices of $G_i$.  Pick two vertices $v\in B, v'\in B'$ of
$X$.  Since $X$ is connected, there are vertices $v_0=v, v_1, \cdots,
v_l=v'$ such that $v_{j-1}$ and $v_j$ are the endpoints of an edge
$e_j$.  Let $B_j$ be the component of $X-X_{T_i}$ that contains $v_j$.
If $e_j\notin E_i$, then $B_{j-1}=B_j$.  On the other hand, if $e_j\in
E_i$, then $v_{B_{j-1}}$ and $v_{B_j}$ are connected by an edge in
$G_i$.

    We next describe the maps from edge spaces to vertex spaces.  Let
 $e_Y$ be an edge of $G_i$ connecting $v_{B_0}$ and $v_{B_1}$. We may
 assume that for each fixed $y\in Y$, $(y,t)$ (where $(y,t)\in Y\times
 (0,1)\subset X$) converges to a point in $B_0$ as $t\ra 0$ and to a
 point in $B_1$ as $t\ra 1$.  The maps $g_{e_Y,0}: Y\ra B_0$ and
 $g_{e_Y,1}: Y\ra B_1$ can be described as follows.  Recall each
 $(n-1)$-cube of $Y$ has the form $C_i$, where $C$ is an $n$-cube of
 $X$ and $C_i\subset C$ is the hyperplane of $C$ containing the
 midpoint of some edge $e\in E_i$ of $C$ and perpendicular to $e$.
 Clearly $C$ has the decomposition $C_i\times [0,1]$. We may assume
 $C_i\times \{0\}\subset C$ is contained in $B_0$ and $C_i\times
 \{1\}\subset C$ is contained in $B_1$.  Then the map $g_{e_Y,0}: Y\ra
 B_0$ sends $C_i$ to $C_i\times \{0\}$ via the identity map.
 Similarly for the map $g_{e_Y,1}$. Thus $g_{e_Y,0}$ and $g_{e_Y,1}$
 are nondegenerate combinatorial maps between FCCs.  Recall that two
 cubes in $X$ either are disjoint or intersect in a single cube.  It
 follows that the maps $g_{e_Y,0}$ and $g_{e_Y,1}$ are immersions,
 that is, they are locally injective.

\begin{Le}\label{localinj}
{The maps $g_{e_Y,0}$ and $g_{e_Y,1}$ are local isometric
 embeddings. In particular, they induce injective homomorphisms
 between fundamental groups.}
\end{Le}
\begin{proof}
We show that $g_{e_Y,0}$ is a local isometric embedding, the proof for
$g_{e_Y,1}$ is similar.  It suffices to show that $g_{e_Y,0}: Y\ra
B_0$ sends geodesics in $Y$ to geodesics in $B_0$.  Recall that $Y$ is
a component of $H_i$ and the component of $X-X_{T_i}$ containing $Y$
is isometric to $Y\times (0,1)$.  $Y$ can be identified with $Y\times
\{\frac{1}{2}\}$. Let $\sigma:I\ra Y$ be a geodesic.  Then for each
$t$, $0<t<1$, the map $\sigma_t: I\ra X$ with
$\sigma_t(s)=(\sigma(s),t)$ for $s\in I$ is a geodesic in $X$.  Since
$X$ has nonpositive curvature, the limit map $\sigma_0:I\ra X$,
$\sigma_0(s):=\lim_{t\ra 0}\sigma_t(s)$ is also a geodesic in $X$(see
p.121, Corollary 7.58 of \cite{BH}).  Now the lemma follows since
$\sigma_0=g_{e_Y,0}\circ \sigma$.
\end{proof}

Summarizing the above observations we have the following:

\begin{Prop}\label{graphdeco}
{Let $X$ be a \FCC of dimension $n$. Then $X$ admits the structure of
 a graph of spaces, where all the vertex and edge spaces are
 $(n-1)$-dimensional \FCCs and the maps from edge spaces to vertex
 spaces are combinatorial immersions.  In particular, $\pi_1(X)$ has a
 graph of groups decomposition.}
\end{Prop}

Let $X$, $G_i$ be as above, $v_B\in G_i$ a vertex and $e_Y\subset G_i$
 an edge incident to $v_B$.  Suppose for each $y\in Y$, $(y,t)$
 converges to a point in $B$ as $t\ra 0$, and let $g_{e_Y,0}: Y\ra B$
 be the map from the edge space to the vertex space.  For each vertex
 $w\in Y$, $\overline{\{w\}\times (0,1)}$ is an edge in $X$; we let
 $e_w$ be the associated oriented edge which goes from $0$ to $1$.
 For any oriented edge $e$ of $X$ with initial point $v$,
 $\stackrel{\ra}{e}\in \Link(X,v)$ denotes the point representing $e$.

\begin{Le}\label{cover}
{In the above notation, the following two conditions are
equivalent:

\emph{(1)}\qua $g_{e_Y,0}$ is not a covering
map;

\emph{(2)}\qua there exists a vertex $w\in Y$ and an oriented
edge $e\subset B$ with initial point $v:=g_{e_Y,0}(w)$ such that
$d_v(\stackrel{\ra}{e}, \stackrel{\ra}{e_w})\ge \pi$.}

\end{Le}

\begin{proof}
We first notice that (2) is equivalent to the following
condition:

(3)\qua There is a vertex $w\in Y$ and an edge
$e\subset B$ incident to $v=g_{e_Y,0}(w)$ such that no edge (in $Y$)
incident to $w$ is mapped to $e$.

So (2) clearly implies (1).
Now assume $g_{e_Y,0}$ is not a covering map.  If the image of
$g_{e_Y,0}$ does not contain all the vertices of $B$, the
connectedness of the $1$-skeleton of $B$ implies there is an edge
$e\subset B$ with two endpoints $v_1$ and $v_2$ such that $v_1\in
\image(g_{e_Y,0})$ and $v_2\notin \image(g_{e_Y,0})$.  In this case no
edge of $Y$ is mapped to $e$ and (3) holds.  So we assume
$\image(g_{e_Y,0})$ contains all the vertices of $B$.  If (3) does not
hold, then $\image(g_{e_Y,0})=B$ and for each vertex $w\in Y$,
$g_{e_Y,0}$ maps $star(w)$ isomorphically onto $star(v)$, where
$v=g_{e_Y,0}(w)$. It follows that $g_{e_Y,0}$ is a covering map,
contradicting to the assumption.
\end{proof}

\subsection{Davis complexes of right-angled Coxeter groups} \label{davis}

In this section we give examples of finite FCCs whose universal covers
 are Davis complexes of certain right-angled Coxeter groups.

Let $S$ be a finite set, and $M=(m_{s,s'})_{s,s'\in S}$ a symmetric
matrix such that $m_{s,s}=1$ for $s\in S$ and $m_{s,s'}\in
\{2,3,\cdots, \}\cup \{\infty\}$ for $s\not=s'\in S$.  The group $W$
given by the presentation $W=<S|(ss')^{m_{s,s'}}=1, s,s'\in S>$ is a
Coxeter group, where $(ss')^\infty=1$ means the relation is void.  The
Coxeter group $W$ is a right-angled Coxeter group if for any
$s\not=s'\in S$, $m_{s,s'}=2$ or $\infty$.

Given any Coxeter group $W$, there is a locally finite cell complex
  $D_W$, the so called Davis complex of $W$ such that $W$ acts on
  $D_W$ properly with compact quotient \cite{D}.  Moussong
  (\cite{M}, \cite{D}) showed that there is a piecewise Euclidean
  metric on $D_W$ that turns $D_W$ into a $\CAT(0)$ space.  The action
  of $W$ on $D_W$ preserves the piecewise Euclidean metric, that is,
  $W$ acts on $D_W$ as a group of isometries.  When $W$ is a
  right-angled Coxeter group, the Davis complex $D_W$ is a cubical
  complex.  Below we shall describe FCCs which are finite quotients of
  Davis complexes of certain right-angled Coxeter groups.

There is a one-to-one correspondence between right-angled Coxeter
     groups and finite flag complexes.  Let $W$ be a right-angled
     Coxeter group with standard generating set $S$.  The nerve $N(W)$
     of $W$ is a simplicial complex with set of vertices $S$.  For any
     nonempty subset $T\subset S$ there is a simplex in $N(W)$ with
     $T$ as its vertex set if and only if $m_{t,t'}=2$ for any
     $t\not=t'\in T$.  $N(W)$ is clearly a flag complex.  Conversely,
     let $K$ be a finite flag complex with set of vertices $S$.  Then
     we can define a right-angled Coxeter group with generating set
     $S$ as follows: for $s\not=s'$, $m_{s,s'}=2$ if there is an edge
     in $K$ joining $s$ and $s'$, and $m_{s,s'}=\infty$ otherwise.

Let $K$ be a finite flag complex with vertex set $S$.  We shall
   construct a finite cubical complex $Y(K)$ whose vertex links are
   all isomorphic to $K$ (\cite{BH}, p.212).  Let $V$ be a Euclidean
   space with dimension equal to $|S|$, the cardinality of $S$.
   Identify the standard basis $e_s$($s\in S$) with $S$ and consider
   the unit cube $[0,1]^{|S|}\subset V$.  For each nonempty subset
   $T\subset S$, let $F_T$ be the face of $[0,1]^{|S|}$ spanned by the
   unit vectors $e_t$, $t\in T$.  The cubical complex $Y(K)$ is the
   subcomplex of $[0,1]^{|S|}$ consisting of all faces parallel to
   $F_T$, for all nonempty subsets $T$ of $S$ that are the vertex sets
   of simplices of $K$.  Notice that $Y(K)$ contains the 1-skeleton of
   the unit cube $[0,1]^{|S|}\subset V$.  In particular, $Y(K)$ is
   connected.

\begin{Prop}\label{davisquo}{\rm\cite{BH}}\qua 
{Let $K$ be a finite flag complex. Then $Y(K)$ is a connected finite
cubical complex with nonpositive curvature all of whose vertex links
are isomorphic to $K$.}
\end{Prop}

We subdivide $Y(K)$ by using the hyperplanes $x_s=1/2$ ($s\in S$) of
the unit cube $[0,1]^{|S|}$.  Let $X(K)$ be the obtained cubical
complex.  Then $X(K)$ is also a finite cubical complex of nonpositive
curvature with some of its vertex links isomorphic to $K$.

A simplicial complex $K$ is foldable if there is a simplicial map
 $f:K\ra \Delta$ onto an $n$-simplex such that the restriction of $f$
 to each simplex is injective.

\begin{Cor}\label{ggg}
{Let $K$ be a finite flag complex. If $K$ is foldable, dimensionally
 homogeneous and has no boundary, then $X(K)$ is a finite \FCC.}
\end{Cor}

\begin{proof}  It is not hard to see that $X(K)$  
 is dimensionally homogeneous and has no boundary.  Proposition
  \ref{gcomplete} implies that $X(K)$ is geodesically complete.  We
  need to show that $X(K)$ is foldable.  First we notice that the
  group $\mathbb{Z}_2^{|S|}$ acts on $X(K)$ as a group of isometries:
  the $s$-th factor $\mathbb{Z}_2$ acts as the orthogonal reflection
  about the hyperplane $x_s=1/2$.  Let $o\in X(K)$ be the origin of
  $V$ and $star(o)$ the star of $o$ in $X(K)$.  Then the quotient of
  $X(K)$ by $\mathbb{Z}_2^{|S|}$ is isomorphic to $star(o)$.  So we
  have a nondegenerate combinatorial map from $X(K)$ onto $star(o)$.
  Since the link of $X(K)$ at $o$ is isomorphic to $K$ and $K$ is
  foldable, the star $star(o)$ can be folded according to the folding
  of $K$.  The composition of these two maps is a folding of $X(K)$.
\end{proof}

Next we construct certain flag complexes that satisfy the assumptions
  in Corollary \ref{ggg}.  Let $n\ge 2$. A standard $n$-sphere is the
  unit round $n$-sphere with an all right simplicial complex
  structure.  A standard $n$-sphere has $n+1$ subcomplexes which are
  standard $(n-1)$-spheres; we call them equators.  Let ${\mathbb
  S}^n$ be a standard $n$-sphere and $E$ one of its equators.  An all
  right simplicial complex is called a standard $n$-hemisphere if it
  is isomorphic to the closure of one of the components of ${\mathbb
  S}^n-E$; the subcomplex of the standard $n$-hemisphere corresponding
  to $E$ is also called its equator.  The unique point on a standard
  $n$-hemisphere that has distance $\pi/2$ from its equator is called
  its pole.

Let ${\mathbb S}^n$ be a standard $n$-sphere and $E_1,\cdots, E_{n+1}$
  its equators.  A \emph{hemispherex} is an all right simplicial
  complex obtained from ${\mathbb S}^n$ by attaching a finite number
  of standard $n$-hemispheres along the equators of ${\mathbb S}^n$
  such that for each $i, 1\le i\le n+1$, there is at least one
  hemisphere attached along $E_i$.  Here the attaching is realized
  through isomorphisms between equators of ${\mathbb S}^n$ and those
  of the hemispheres.  It is clear that a hemispherex $H$ satisfies
  all the conditions in Corollary \ref{ggg}, so the corresponding
  $X(H)$ is a finite FCC.  The universal cover of $X(H)$ is a
  subdivision of the Davis complex of the right-angled Coxeter group
  whose nerve is $H$.

 Hemispherex was first introduced by Ballmann and Brin in \cite{BBr2}.
 The Euclidean cone over a hemispherex is the first example of a
 higher rank $\CAT(0)$ space aside from Euclidean buildings and
 symmetric spaces.  But such a space does not admit cocompact
 isometric actions. On the other hand, we shall see in Section
 \ref{rankonegeo} that if $X$ is a finite FCC with some vertex link
 isomorphic to a hemispherex, then it contains closed rank one
 geodesics. In particular, $X(H)$ contains closed rank one geodesics
 if $H$ is a hemispherex.

\section{Existence of closed rank one geodesics}\label{existence}

In this section we discuss the existence of closed rank one geodesics
in a finite FCC.  Throughout this section $X$ denotes a finite FCC of
dimension $n$, except in Section \ref{rankrigidity}.

\subsection{Vertex links and  $1$-skeleton} \label{vlink}

Recall that $X$ is dimensionally homogeneous, geodesically complete
and has nonpositive curvature.  It follows that for each vertex $v$ of
$X$, the link $\Link(X,v)$ is dimensionally homogeneous, has no
boundary and is a flag complex.

Recall that the set of edges in $X$ is a disjoint union:
 $E=\amalg_{i=1}^n E_i$. Let $v$ be a vertex of $X$. For $1\le i\le
 n$, let $V_{i,v}$ be the set of vertices in $\Link(X,v)$ that
 correspond to edges in $E_i$.  Since $X$ has no boundary, $V_{i,v}$
 contains at least two points.

\begin{Le}\label{twolink}
{Let $v\in X$ be a vertex and $1\le i\le n$. If $V_{i,v}$ consists of
 exactly two points, then $\Link(X,v)$ is the spherical join of
 $V_{i,v}$ and $L_{i,v}$, where $L_{i,v}$ is the subcomplex of
 $\Link(X,v)$ consisting of all simplices of $\Link(X,v)$ with vertices
 in $V_{j,v}$, $j\not=i$.}
\end{Le}

\begin{proof}
Denote by $\xi_+, \xi_- $ the two vertices in $V_{i,v}$, and let
  $\eta\in V_{j,v}$ with $j\not=i$.  Since $\Link(X,v)$ is a flag
  complex it suffices to show that there are edges $\eta\xi_+$ and
  $\eta\xi_-$ in $\Link(X,v)$.

Since $\Link(X,v)$ is dimensionally homogeneous, there is a
 $(n-1)$-simplex $\Delta_1$ of $\Link(X,v)$ containing $\eta$.
 $\Delta_1$ contains exactly one vertex from $V_{i,v}$. Without loss
 of generality we may assume it is $\xi_+$ (thus there is an edge
 $\eta\xi_+$).  Let $\Delta_1'\subset \Delta_1$ be the $(n-2)$-face
 disjoint from $\xi_+$.  Since $\Link(X,v)$ has no boundary, there is a
 $(n-1)$-simplex $\Delta_2\not=\Delta_1$ containing $\Delta_1'$.  The
 assumption that $V_{i,v}$ consists of exactly two points implies
 $\xi_-$ is a vertex of $\Delta_2$.  In particular there is an edge
 $\eta\xi_-$.
\end{proof}

 For any vertex $v$ of $X$ and $1\le i\le n$, let $X_{i,v}$ be the
 component of $X_{\{i\}}$ that contains $v$.

\begin{Cor}\label{piangle}
{Let $v\in X$ be a vertex and $1\le i\le n$.  If there are $\xi\in
 V_{j,v}$, $\eta\in V_{i,v}$ with $j\not=i$ and $d_v(\xi, \eta)\ge
 \pi$, then $X_{i, v}$ is not a circle.

}
\end{Cor}

\begin{proof}
By Lemma \ref{twolink}, $V_{i,v}$ contains at least three points and
  therefore there are at least three edges of $X_{i, v}$ incident to
  $v$.
\end{proof}

For an oriented edge $e$, $\bar e$ denotes the same edge with the
 opposite orientation and $t(e)$ denotes the terminal point of $e$.
 It is not hard to prove the following lemma.

\begin{Le}\label{graph}
{Let $\Gamma$ be a connected finite graph such that the valence of
  each vertex is at least two.  Assume $\Gamma$ is not homeomorphic to
  a circle.  Then for any two oriented edges $e_1$, $e_2$ in $\Gamma$,
  there is a geodesic $c$ from $t(e_1)$ to $t(e_2)$ such that $\bar
  e_2*c*e_1$ is also a geodesic.  In particular, for any oriented edge
  $e$ of $\Gamma$, there is a geodesic loop $c$ based at $t(e)$ such
  that $\bar e*c*e$ is also a geodesic.}
\end{Le}

\subsection{Closed rank one geodesics} \label{rankonegeo}

Let $v\in X$ be a vertex and $\xi\in \Link(X,v)$.  Then $\xi$ lies in
  the interior of a simplex $\Delta$ of $\Link(X,v)$.  Define
  $T(\xi)=\{i: V_{i,v}\cap \Delta\not=\phi\}$.  Then $T(\xi)\subset
  \{1,2, \cdots,n\}$. We say $T(\xi)$ is the type of $\xi$.  For an
  oriented edge $e$ of $X$ with initial point $v$, define
  $P_e=\{\xi\in \Link(X,v): d_v(\xi, \stackrel{\ra}{e})=\pi/2\}$.

\begin{Le}\label{transfer}
{Let $e$ be an oriented edge of $X$. Then there is an isometry $D_e:
  P_e\ra P_{\bar e}$ such that $\xi$ and $D_e(\xi)$ have the same type
  for any $\xi\in P_e$. }
\end{Le}

\begin{proof}
Let $i$ be the index such that the geometric edge of $e$ lies in $
 E_i$, $m$ the midpoint of $e$ and $Y$ the component of $H_i$
 containing $m$.  Denote by $v$ and $w$ the initial and terminal
 points of $e$, and $B_0$ and $B_1$ the components of $X_{T_i}$
 containing $v$ and $w$ respectively.  Since $g_{e_Y,0}: Y\ra B_0$ and
 $g_{e_Y,1}: Y\ra B_1$ are local isometric embeddings and $B_0$, $B_1$
 are locally convex in $X$, the induced maps $h_0:\Link(Y,m)\ra
 \Link(B_0,v)\subset \Link(X,v)$ and $h_1: \Link(Y,m)\ra \Link(B_1,
 w)\subset \Link(X,w)$ are isometric embeddings.  It is not hard to
 check that the images of these maps are $P_e$ and $P_{\bar e}$
 respectively.  Set $D_e=h_1\circ {h_0}^{-1}$.  Then $D_e: P_e\ra
 P_{\bar e}$ is an isometry.  Since $h_0$ and $h_1$ clearly preserve
 type, $D_e$ also preserves type.
\end{proof}

\begin{Cor}\label{ttt}
{Let $X$ be a \FCC, $v\in X$ a vertex and $e_1,e_2, e_3$ three
 oriented edges with initial point $v$.  Suppose
 $d_v(\stackrel{\ra}{e_1},
 \stackrel{\ra}{e_3})=d_v(\stackrel{\ra}{e_2}, \stackrel{\ra}{e_3})=
 \pi/2$.  Then \linebreak $d_v(\stackrel{\ra}{e_1},
 \stackrel{\ra}{e_2})=\pi$ if and only if
 $d_w(D_{e_3}(\stackrel{\ra}{e_1}),
 D_{e_3}(\stackrel{\ra}{e_2}))=\pi$, where $w$ is the terminal point
 of $e_3$.}
\end{Cor}

\begin{proof}
It follows easily from the fact that $D_{e_3}:P_{e_3} \ra P_{\bar e_3}
   $ is an isometry.
\end{proof}

We say a geodesic in a metric space with nonpositive curvature has
 rank one if its lifts in the universal cover have rank one.

\begin{Prop}\label{ncolor}
{Let $c\subset X$ be a closed geodesic that is contained in the
 $1$-skeleton. If for each $i$, $c$ contains at least one edge from
 $E_i$, then $c$ is a closed rank one geodesic.}
\end{Prop}
\begin{proof} 
Let $\pi:\tilde X\ra X$ be the universal cover of $X$, $\tilde c$ a
   lift of $c$ to $\tilde X$ and $\tilde E_i=\pi^{-1}(E_i)$ for $1\le
   i\le n$.  Note $\tilde X$ is also a FCC and the set $\tilde E$ of
   edges of $\tilde X$ has the decomposition into different colors
   $\tilde E=\amalg_{i=1}^{n}\tilde E_i$.

Assume the proposition is false.  Then there is half flatplane $H$
bounding $\tilde c$, that is, $H=\image(f)$ where $f$ is an isometric
embedding $f:\R\times [0,\i)\ra \tilde X$ with $f(\R\times
\{0\})=\tilde c$.  For any vertex $v\in \tilde c$, $H$ determines a
unique point $\xi_v\in \Link(\tilde X, v)$ with the following property:
for any oriented edge $e\subset \tilde c$ with initial point $v$ and
terminal point $w$, $\xi_v\in P_e$ and $D_e(\xi_v)=\xi_w$.  By Lemma
\ref{transfer}, all $\xi_v$ have the same type $T\subset \{1,2,\cdots,
n\}$.  Let $i\in T$. By assumption there is an oriented edge $e\subset
\tilde c$ with $e\in \tilde E_i$.  Let $v$ be the initial point of
$e$.  Since $\Link(\tilde X, v)$ is an all right spherical complex and
$\stackrel{\ra}{e}\in \Link(\tilde X, v)$ is a vertex, $\xi_v\in P_e$
implies $d_v(\stackrel{\ra}{e}, \xi)=\pi/2$ for any vertex $\xi$ of
the smallest simplex of $\Link(\tilde X, v)$ containing $\xi_v$.  Since
$T(\xi_v)=T\owns i$, the definition of type implies there is an
oriented edge $e_i\in \tilde E_i$ with $d_v(\stackrel{\ra}{e},
\stackrel{\ra}{e_i})=\pi/2$.  This contradicts to the facts that $e\in
\tilde E_i$ and $\tilde X_{i,v}$ is locally convex in $\tilde X$.
\end{proof}

Let $v\in X$ be a vertex.  We define a relation $\sim_v$ on the set
 $\{1,2,\cdots, n\}$ as follows: $i\sim_v j$ if and only if there are
 oriented edges $e_1\in E_i$, $e_2\in E_j$ with initial point $v$ such
 that $d_v(\stackrel{\ra}{e_1}, \stackrel{\ra}{e_2})\ge \pi$.

\begin{Prop}\label{equi}
{Let $X$ be a finite \FCC of dimension $n$.  For a fixed vertex $v\in
   X$, if $T\subset \{1,2,\cdots, n\}$ is an equivalence class with
   respect to the equivalence relation generated by $\sim_v$, then
   there is a closed geodesic in the $1$-skeleton of $X$ that contains
   at least one edge from $E_i$ for each $i\in T$.}
\end{Prop}

\begin{proof}
By the definition of the equivalence relation, we may assume $T=\{i_1,
   i_2, \cdots, i_m\}$ such that for each $t$, $2\le t\le m$, there is
   some $j$, $1\le j<t$ with $i_t\sim_v i_j$.  Now we prove the
   following claim by induction on $k$, $1\le k\le m$: there is a
   closed geodesic $c_k$ in the $1$-skeleton of $X$ such that for each
   $j\le k$, $c_k$ contains an edge belonging to $E_{i_j}$ and
   incident to $v$.

The claim is clear for $k=1$. Now let $k\ge 2$ and assume the claim
 has been established for $k-1$.  By the above paragraph we have
 $i_k\sim_v i_j$ for some $j<k$.  Hence there are oriented edges $
 e_j'\in E_{i_j}$, $e_k\in E_{i_k}$ with initial point $v$ such that
 $d_v(\stackrel{\ra}{e_j'}, \stackrel{\ra}{e_k})\ge \pi$.  Corollary
 \ref{piangle} implies $X_{i_j,v}$ and $X_{i_k,v}$ are not circles.
 By induction hypothesis, there is a closed geodesic $c_{k-1}$ in the
 $1$-skeleton of $X$ containing an oriented edge $e_j\in E_{i_j}$ with
 initial point $v$.  Let $v_1$, $v_2$, $v_3$ be the terminal points of
 the edges $e_k$, $ e_j'$, $e_j$ respectively.  By Lemma \ref{graph},
 there is a geodesic loop $\tilde c_1$ in $X_{i_k,v}$ based at $v_1$
 such that $c_1'=\bar e_k*\tilde c_1*e_k$ is also a geodesic in
 $X_{i_k,v}$. Similarly there are geodesic loops $\tilde c_2$ and
 $\tilde c_3$ in $X_{i_j,v}$ based at $v_2$ and $v_3$ respectively
 such that $c_2'=\bar e_j'*\tilde c_2*e_j'$ and $c_3'=\bar e_j*\tilde
 c_3* e_j$ are geodesics in $X_{i_j,v}$.  Proposition \ref{convex}
 implies $c_1'$, $c_2'$ and $c_3'$ are geodesics in $X$.  We
 reparametrize the closed geodesic $c_{k-1}$ so that it starts from
 $v$ with initial segment $e_j$.  Define
 $c_k=c_2'*c_3'*c_{k-1}*c_2'*c_1' $ if $e_j'\not=e_j$ and
 $c_k=c_1'*c_2'*c_{k-1}$ if $e_j'=e_j$. Now it is easy to check that
 $c_k$ is a closed geodesic with the required property.
\end{proof}

The following corollary follows immediately from Propositions
\ref{ncolor} and \ref{equi}.

\begin{Cor}\label{nco}
{Let $X$ be a finite \FCC of dimension $n$ and $v\in X$ a vertex. If $
\{1,2,\cdots, n\}$ is a single equivalence class with respect to the
equivalence relation generated by $\sim_v$, then there is a closed
rank one geodesic contained in the $1$-skeleton of $X$.}

\end{Cor}

Let $H$ be a hemispherex with central sphere ${\mathbb S}^n$, and
 $E_1, \cdots, E_{n+1}$ the equators of ${\mathbb S}^n$. For each $i$,
 $1\le i\le n+1$, let $H_i$ be a fixed hemisphere of $H$ that is
 attached to ${\mathbb S}^n$ along $E_i$. Denote the pole of $H_i$ by
 $p_i$.  Then the distance between $p_i$ and $p_j$ is $\pi$ for
 $i\not=j$.  Now let $X$ be a FCC and $v\in X$ a vertex with
 $\Link(X,v)=H$.  Then there are oriented edges $e_i$, $1\le i\le n+1$
 with initial point $v$ and $\stackrel{\ra}{e_i}=p_i$.  It is clear
 that $e_i$ and $e_j$ have different colors for $i\not=j$. It follows
 that the assumption in the above corollary is satisfied if the vertex
 link $\Link(X,v)$ is a hemispherex. Thus we have:

\begin{Cor}\label{heme}
{Let $X$ be a finite \FCC of dimension $n$. Suppose there is a vertex
$v\in X$ such that $\Link(X,v)$ is a hemispherex, then $X$ has a closed
rank one geodesic contained in the $1$-skeleton.}
\end{Cor}

\begin{Prop}\label{largepi}
{Let $X$ be a finite \FCC. Suppose there is a vertex $v\in X$ and
 $\xi\in V_{i,v}$, $\eta\in V_{j,v}$ with $i\not=j$ and $d_v(\xi,
 \eta)>\pi$, then $X$ contains a closed rank one geodesic in the
 $1$-skeleton. }

\end{Prop}

\begin{proof}
Let $e_1$ and $e_2$ be the two oriented edges with initial point $v$
 that give rise to $\xi$ and $\eta$ respectively.  Since
 $d_v(\xi,\eta)>\pi$, Corollary \ref{piangle} implies $X_{i,v}$ and
 $X_{j,v}$ are not circles.  Lemma \ref{graph} then implies there are
 geodesic loops $c_1\subset X_{i,v}$ and $c_2\subset X_{j,v}$ based at
 $t(e_1)$ and $t(e_2)$ respectively such that $c_1':=\bar e_1*c_1*e_1$
 and $c_2':=\bar e_2*c_2*e_2$ are geodesics in $X_{i,v}$ and $X_{j,v}$
 respectively.  Since by Proposition \ref{convex} $X_{i,v}$ and
 $X_{j,v}$ are locally convex in $X$, $c_1'$ and $c_2'$ are geodesics
 in $X$.  Let $c= c_2'*c_1'$.  Since $d_v(\xi,\eta)>\pi$, it is clear
 that $c$ is a closed rank one geodesic.
\end{proof}

\subsection{A splitting criterion} \label{splitting}

Let $Y$ be a $\CAT(0)$ space and $Z_1, Z_2\subset Y$ be two closed,
 convex subsets.  We say $Z_1$, $Z_2$ are \emph{parallel} if for some
 $a\ge 0$ there is an isometric embedding $f: Z_1\times [0,a]\ra Y$
 such that $f(Z_1\times \{0\})=Z_1$ and $f(Z_1\times \{a\})=Z_2$.  For
 any closed convex subset $Z\subset Y$ of a $\CAT(0)$ space $Y$, let
 $P_Z$ be the union of all closed convex subsets that are parallel to
 $Z$.  When $Z$ is geodesically complete, $P_Z$ is closed, convex and
 isometrically splits $Z\times C$, where $C\subset Y$ is closed and
 convex (\cite{BBr2}, p.6).

\begin{Prop}\label{split}
{Let $X$ be a \FCC of dimension $n$.  Suppose $\{1,2,\cdots, n\}$ is
the disjoint union of nonempty subsets $T$, $S$ with the following
property: for any vertex $v\in X$, and any two edges incident to $v$,
$e_i\in E_i$, $e_j\in E_j$ with $i\in T$, $j\in S$, there is a square
containing $e_i$ and $e_j$ in the boundary.  Then the universal cover
$\tilde X$ of $X$ is isometric to the product of two $\CAT(0)$ \FCCs.
}

\end{Prop}

\begin{proof}
For any vertex $v\in \tilde X$, let $\tilde X_{T,v}$ be the component
  of $\tilde X_T$ that contains $v$.  We claim for any edge $e$ of
  $\tilde X$ with endpoints $v$ and $w$, $\tilde X_{T,v}$ and $\tilde
  X_{T,w}$ are parallel.

 We may assume $e\in E_i$ for some $i\in S$, otherwise $\tilde
 X_{T,v}=\tilde X_{T,w}$.  For any $k\ge 0$ we inductively define a
 subcomplex $\tilde X_{T,v}(k)$ of $ \tilde X_{T,v}$: $\tilde
 X_{T,v}(0)=\{v\}$, for $k\ge 1$, $\tilde X_{T,v}(k)$ is the union of
 $\tilde X_{T,v}(k-1)$ and all the cubes in $\tilde X_{T,v}$ that have
 nonempty intersection with $\tilde X_{T,v}(k-1)$.  Similarly one can
 define $\tilde X_{T,w}(k)$.  We also define subcomplexes $\tilde
 X_{T,e}(k)$ of $\tilde X$: $\tilde X_{T,e}(0)=e$, for $k\ge 1$,
 $\tilde X_{T,e}(k)$ is the union of $\tilde X_{T,e}(k-1)$ and all the
 cubes whose edges are in $E_i \cup(\cup_{j\in T}E_j)$ and whose
 intersections with $\tilde X_{T,e}(k-1)$ contain edges from $E_i$.
 Set $\tilde X_{T,e}=\cup_{k\ge 0}\tilde X_{T,e}(k)$.

Since $\tilde X$ is a $\CAT(0)$ cubical complex, the vertex links of
   $\tilde X$ are flag complexes.  Our assumption then implies that
   for any $(m-1)$ ($m\le n$) cube $C$ in $\tilde X_{T,v}(1)$ or
   $\tilde X_{T,w}(1)$, there is a unique $m$-cube in $\tilde X$ that
   contains both $e$ and $C$.  It follows that $\tilde X_{T,e}(1)$
   contains $\tilde X_{T,v}(1)$ and $\tilde X_{T,w}(1)$ and there is
   an isomorphism
 $$f_{e,1}: \tilde X_{T,v}(1)\times [0,1] \ra \tilde X_{T,e}(1) $$
such that ${f_{e,1}}_{|\tilde X_{T,v}(1)\times \{0\}}$ is the identity
map and $f_{e,1}(\tilde X_{T,v}(1)\times \{1\})=\tilde X_{T,w}(1)$.
Now $\tilde X_{T,v}(k)=\cup_{v'}\tilde X_{T,v'}(1)$ and $\tilde
X_{T,e}(k)=\cup_{e'}\tilde X_{T,e'}(1)$, where $v'$ varies over all
vertices in $\tilde X_{T,v}(k-1)$ and $e'\subset \tilde X_{T,e}(k-1)$
varies over all edges from $E_i$.  Notice all the maps $f_{e',1}$ are
compatible for $e'\subset \tilde X_{T,e}(k-1)$ from $E_i$.  It follows
that for each $k$ there is an isomorphism $f_{e,k}:\tilde
X_{T,v}(k)\times [0,1] \ra \tilde X_{T,e}(k) $ such that
${f_{e,k}}_{|\tilde X_{T,v}(k)\times \{0\}}$ is the identity map,
$f_{e,k}(\tilde X_{T,v}(k)\times \{1\})=\tilde X_{T,w}(k)$ and
$f_{e,k}$ agrees with $f_{e,k-1}$ when restricted to $\tilde
X_{T,v}(k-1)\times [0,1]$.  The union of all these isomorphisms
$f_{e,k}$ defines an isomorphism $f_e:\tilde X_{T,v}\times [0,1] \ra
\tilde X_{T,e} $ such that ${f_e}_{|\tilde X_{T,v}\times \{0\}}$ is
the identity map onto $\tilde X_{T,v}$ and ${f_e}({\tilde
X_{T,v}\times \{1\}})=\tilde X_{T,w}$.  It follows that $\tilde
X_{T,v}$ and $\tilde X_{T,w}$ are parallel.

Fix a vertex $v_0\in \tilde X$ and let $P_T$ be the parallel set of
$\tilde X_{T,v_0}$.  Note $\tilde X_{T,v_0}$ is closed, convex and
geodesically complete.  It follows that $P_T$ isometrically splits
$P_T=\tilde X_{T,v_0}\times Y$ where $Y\subset \tilde X$ is a closed
convex subset and for $y\in Y$, $\tilde X_{T,v_0}\times \{y\}$ is
parallel to $\tilde X_{T,v_0}$.  By the claim we have established, all
vertices of $\tilde X$ lie in $P_T$. Since $\tilde X$ is the convex
hull of all its vertices we see $\tilde X=P_T$ splits.  \end{proof}

Proposition \ref{split} has previously been established in dimension 2
(Theorem 1.10 on p.36 of \cite{W} and Theorem 10.2 in \cite{BW}).

Let $X$ be a FCC.  By Proposition \ref{graphdeco} $X$ has the
   structure of a graph of spaces, where all the vertex and edge
   spaces are FCCs and the maps from edge spaces to vertex spaces are
   combinatorial immersions.  The following corollary follows from
   Lemma \ref{cover} and Proposition \ref{split}.

\begin{Cor}\label{covering}
{Let $X$ be a \FCC with a graph of spaces decomposition as in
  Proposition \ref{graphdeco}.  If all the maps from edge spaces to
  vertex spaces are covering maps, then the universal cover of $X$ is
  isometric to the product of a simplicial tree and a
  $(n-1)$-dimensional $\CAT(0)$ \FCC.}
\end{Cor}

\subsection{Rank rigidity in low dimensions} \label{rankrigidity}

In this section we discuss the rank rigidity problem for finite FCCs
with dimension $\le 3$.  A 1-dimensional finite FCC $X$ is a finite
graph and each of its vertices is incident to at least two edges; $X$
clearly contains closed geodesics and all the geodesics in $X$ have
rank one.  The claim in dimension 2 follows easily from Corollary
\ref{nco} and Proposition \ref{split}.

\begin{Th}\label{rigidity}
{Let $X$ be a finite \FCC of dimension $3$ with universal cover
$\tilde X$.  

\emph{(1)}\qua If $\tilde X$ has higher rank, then $\tilde
X$ is isometric to the product of two $\CAT(0)$ \FCCs. 

\emph{(2)}\qua If
$\tilde X$ has rank one, then there is a closed rank one geodesic in
the $1$-skeleton of $X$.}

\end{Th}

\begin{proof} Suppose  that $\tilde X$ does not split  as  a  product. 
Recall the decomposition of the set of edges into different colors:
  $E=E_1\amalg E_2\amalg E_3$. We shall call the edges in $E_1$,
  $E_2$, $E_3$ blue, green and red edges respectively.  By Proposition
  \ref{largepi} and Proposition \ref{split} we may assume the
  following: for any two oriented edges $e_1$, $e_2$ with the same
  initial point $v$ but different colors, $d_v(\stackrel{\ra}{e_1},
  \stackrel{\ra}{e_2})\le \pi$ holds; there exist two oriented edges
  $e_1$, $e_2$ with different colors (say blue and green respectively)
  and the same initial point $v$ such that $d_v(\stackrel{\ra}{e_1},
  \stackrel{\ra}{e_2})=\pi$.

   Consider the graph of spaces decomposition of $X$ where the vertex
  spaces are components of $X_{\{1,2\}}$.  Let $G_3$ be the base
  graph.  By Corollary \ref{covering} at least one of the maps from
  the edge spaces to vertex spaces is not a covering map.  Recall the
  base graph $G_3$ is connected.  Let $k\ge 0$ be the smallest integer
  with the following property: there are two vertices $v_B$, $v_{B'}$
  of $G_3$ at distance $k$ apart such that 

(1)\qua $B$ contains
  an oriented blue edge $e_1$ and an oriented green edge $e_2$ with
  the same initial point $v\in X$ such that $d_v(\stackrel{\ra}{e_1},
  \stackrel{\ra}{e_2})=\pi$; 

 (2)\qua there is an edge $e_Y\subset
  G_3$ incident to $v_{B'}$ such that the map from the edge space $Y$
  to the vertex space $B'$ is not a covering map.

 We claim $k=0$. Assume $k\ge 1$. Let $B_0=B, B_1, \cdots, B_k=B'$ be
 a sequence of components of $X_{\{1,2\}}$ such that $v_{B_i}$ and
 $v_{B_{i+1}}$ ($0\le i\le k-1$) are adjacent vertices in $G_3$.
 Since $k\ge 1$, the map from the edge space of $v_Bv_{B_1}$ to the
 vertex space $B$ is a covering map.  Lemma \ref{cover} implies that
 there is a red edge $e$ with one endpoint $v$ and the other endpoint
 $w$ in $B_1$ such that $e$ is perpendicular to both $e_1$ and $e_2$.
 Corollary \ref{ttt} implies that there exist an oriented blue edge
 $e_3$ and an oriented green edge $e_4$ with initial point $w$ such
 that $d_w(\stackrel{\ra}{e_3}, \stackrel{\ra}{e_4})=\pi$. Note $e_3,
 e_4\subset B_1$ and the distance from $v_{B_1}$ to $v_{B'}$ in $G_3$
 is $k-1$, contradicting to the definition of $k$.  Therefore $k=0$
 and $B=B'$.  By Lemma \ref{cover} there is a vertex $v'\in B$ and
 oriented edges $e_r$ (red), $e_b\subset B$ with initial point $v'$
 such that $d_{v'}(\stackrel{\ra}{e_r}, \stackrel{\ra}{e_b})=\pi$.  We
 may assume $e_b$ is a blue edge.

By Corollary \ref{nco} we may assume the following: for any vertex
$v\in B$, if there are oriented blue edge $e_1$ and red edge $e_2$
with initial point $v$ such that\linebreak $d_v(\stackrel{\ra}{e_1},
\stackrel{\ra}{e_2})=\pi$, then all green edges incident to $v$ are
perpendicular to all blue and red edges incident to $v$.

  Recall that $B$ is a finite FCC of dimension $2$.  We consider the
graph of spaces decomposition of $B$ where the vertex spaces consist
of blue edges.  Let $G$ be the connected base graph.  Lemma
\ref{cover} and the condition (1) above imply that not all maps from
edge spaces to vertex spaces are covering maps.  Let $l\ge 0$ be the
smallest integer with the following property: there are two vertices
$v_C$ and $v_{C'}$ of $G$ at distance $l$ apart such that

(1)\qua
there is a vertex $v'\in C$, an oriented red edge $e_r$ and an
oriented blue edge $e_b$ with initial point $v'$ such that
$d_{v'}(\stackrel{\ra}{e_r}, \stackrel{\ra}{e_b})=\pi$;

(2)\qua
there is an edge $e_Y\subset G$ incident to $v_{C'}$ such that the map
from the edge space $Y$ to the vertex space $C'$ is not a covering
map.

Now the preceding paragraph and a similar argument as
above show that $l=0$ and $C=C'$.

There is a vertex $v''\in C$, an oriented blue edge $e_{b'}$ and an
  oriented green edge $e_g$ with initial point $v''$ such that
  $d_{v''}(\stackrel{\ra}{e_{b'}}, \stackrel{\ra}{e_g})=\pi$.
  Corollary \ref{piangle} implies $C=X_{1,v''}$ is not a circle. By
  Lemma \ref{graph} there is a geodesic $c\subset C$ from $v'$ to
  $v''$ which starts with $e_b$ and ends with $\bar e_{b'}$.
  Similarly there are geodesic loops $c_1$ and $c_2$ based at $t(e_r)$
  and $t(e_g)$ respectively such that $c_1'=\bar e_r*c_1*e_r$ and
  $c_2'=\bar e_g*c_2*e_g$ are also geodesics.  Set $c'= c_1'*\bar
  c*c_2'*c$. Then $c'$ is a closed geodesic that contains blue, green
  and red edges. By Proposition \ref{ncolor} $c'$ is a rank one
  geodesic.
\end{proof}

Theorem \ref{rigidity} (1) also follows from a theorem of Ballmann and
Brin \cite{BBr2}: they proved that if a $3$-dimensional,
geodesically complete and piecewise Euclidean polyhedra $Y$ is
$\CAT(0)$, has higher rank and admits a cocompact and properly
discontinuous group of isometries, then $Y$ either isometrically
splits or is a thick Euclidean building of type $\tilde {A_3}$ or
$\tilde {B_3}$.  A FCC certainly can not be an Euclidean building of
type $\tilde {A_3}$ or $\tilde {B_3}$.  The main point of Theorem
\ref{rigidity} is the existence of closed rank one geodesics in rank
one finite FCCs.  Our proof of Theorem \ref{rigidity} is independent
of the proof in \cite{BBr2}.

\section{Tits alternative for foldable cubical complexes}\label{titsal}

In this section we give a short proof of the Tits alternative for the
fundamental group of a finite FCC.  Ballmann and Swiatkowski have a
slightly more general result \cite{BSw}.

\begin{Th}\label{tits}
{Let $X$ be a finite \FCC. Then any subgroup of $\pi_1(X)$ either
contains a free group of rank two or is virtually free abelian.  }
\end{Th}

\begin{proof}
 We induct on the dimension of $X$. If $\dim X=1$, then $X$ is a
 finite graph and the Theorem clearly holds.  Let $n=\dim X$ and $H$ a
 subgroup of $\pi_1(X)$.  By Proposition \ref{graphdeco} $X$ admits a
 graph of spaces decomposition where all the vertex and edge spaces
 are $(n-1)$-dimensional finite FCCs. It follows that $\pi_1(X)$
 admits a graph of groups decomposition and acts on the associated
 Bass-Serre tree $T$. As a subgroup of $\pi_1(X)$, $H$ also acts on
 $T$.  By \cite{PV}, $H$ contains a free group of rank two unless one
 of the following happens: $H$ fixes a point in $T$, $H$ stabilizes a
 complete geodesic in $T$, or $H$ fixes a point in $\partial_\infty
 T$.  We need to consider these three exceptional cases.

First assume $H$ fixes a point in $T$.  Then $H$ fixes a vertex of $T$
  and is a subgroup of a conjugate of some vertex group of the graph
  of groups decomposition for $\pi_1(X)$.  We have observed that such
  a vertex group is the fundamental group of a finite FCC with
  dimension $n-1$. By induction hypothesis the claim on $H$ holds.

Now assume $H$ stabilizes a complete geodesic $c$ in $T$. Since $T$ is
  a simplicial tree, by taking an index two subgroup we may assume $H$
  acts on $c$ as translations and so there is an exact sequence $1\ra
  N\ra H\ra \Z\ra 1$.  Thus $N\subset H$ has a fixed point in
  $c\subset T$. By the previous paragraph $N$ contains a free group of
  rank two or is virtually free abelian. We may assume $N$ is
  virtually free abelian.  It implies that $H$ is virtually solvable.
  The claim on $H$ follows since any virtually solvable subgroup of a
  group acting properly and cocompactly by isometries on a $\CAT(0)$
  space is virtually free abelian (\cite{BH}, p.249).  The case when
  $H$ fixes a point in $\partial_\infty T$ can be handled similarly.
\end{proof}

\bibliographystyle{gtart}

\Addresses\recd
\end{document}